\documentclass{article}
\usepackage{amsfonts}
\usepackage{amssymb}
\usepackage{amsmath}
\numberwithin{equation}{section}
\usepackage[mathscr]{eucal}
\usepackage{eufrak}

\begin{document}
\begin{center}
\Large{\textbf{Compact Inverses of The Multipoint Normal Differential Operators For First Order}}\\
\small{Erdal UNLUYOL $^{1*}$ , Elif OTKUN \c{C}EVIK $^{1}$, Zameddin  ISMAILOV }
\footnote [1]{$^{*}$ corresponding author: eunluyol@yahoo.com (Erdal UNLUYOL), Fax: +90 (462) 325 31 95}
\\
\footnotesize{$^{1}$ Karadeniz Technical University, Faculty of Science, Department of Mathematics, 61080 Trabzon, Turkey}
\[ eunluyol@yahoo.com, e_{-}otkuncevik@hotmail.com, zameddin@yahoo.com \]
\end{center}
\centerline{}

\begin{abstract}
In this work, firstly all normal extensions of a multipoint minimal operator generated by linear multipoint differential-operator expression for first order in the Hilbert space of vector functions in terms of boundary values at the endpoints of the infinitely many separated subintervals are described. Finally, a compactness properties of the inverses of such extensions has been investigated.
\end{abstract}

\begin{flushleft}
{\footnotesize \textbf{2010 AMS Subject Classification:} 47A05, 47A20
 }
\end{flushleft}

\section{\textbf{Introduction}}

 It is known that traditional infinite direct sum of Hilbert spaces $H_{n} ,{\rm \; }n\ge 1$ and infinite direct sum of operators $A_{n} $ in $H_{n} ,{\rm \; }n\ge 1$ are defined as
\[ H=\mathop{\oplus }\limits_{n=1}^{\infty } H_{n} =\left\{u=\left(u_{n} \right):u_{n} \in H_{n} ,{\rm \; }n\ge 1{\rm \; and\; }\sum _{n=1}^{\infty }\left\| u_{n} \right\| _{H_{n} }^{2} <+\infty  \right\}, \]
\[A=\mathop{\oplus }\limits_{n=1}^{\infty } A_{n}, D(A)=\left\{ u=(u_{n})\in H:u_{n} \in D(A_{n}),{\rm \; }n\ge 1 {\rm \; and}{\rm \;} Au=\left(A_{n} u_{n} \right)\in H \right\}. \]
\indent A linear space $H$ is a Hilbert space with norm induced by the inner product
\[\left(u,v\right)_{H}=\sum_{n=1}^{\infty}\left(u_{n},v_{n}\right)_{H_{n}},{\rm \;}u, v\in H{\rm \;} \cite{Dun}. \]
\indent The general theory of linear closed operators in Hilbert spaces and its applications to physical problems have been investigated by many researches (for example, see \cite{Dun}, \cite{RoKh}).\\
\indent However, many physical problems requires the study of the theory of linear operators in direct sums in Hilbert spaces ( for example, see \cite{Tim}-\cite{So} and references therein).\\
\indent We note that a detail analysis of normal subspaces and operators in Hilbert spaces have been studied in \cite{Cod} (see references in it ).\\
\indent This study contains three section except introduction. In section 2, the multipoint minimal and maximal operators for the first order differential-operator expression are determined. In section 3, all normal extensions of multipoint formally normal operators are described in terms of boundary values in the endpoints of the infinitely many separated subintervals. Finally In section 4, compactness properties of the inverses of such extensions have been established.

\section{\textbf{The Minimal and Maximal Operators}}

 Throughout this work $\left(a_{n} \right)$ and $\left(b_{n} \right)$ will be sequences of real numbers such that
\[-\infty <a_{n} <b_{n} < a_{n+1} <\cdots <+\infty ,\]
\noindent $H_{n} $ is any Hilbert space, {\rm \;\;\;}$\Delta _{n} =\left(a_{n} ,b_{n} \right), {\rm \;\;\;}L_{n}^{2} =L^{2} \left(H_{n}, {\rm \;\;\;}\Delta _{n} \right), L^{2} =\mathop{\oplus }\limits_{n=1}^{\infty } L^{2}\left(H_{n} ,\Delta _{n} \right),{\rm \;\;\;} n\ge 1,\\
\noindent \underset{n\geq1}{sup}(b_{n}-a_{n})<+\infty, W_{2}^{1} =\mathop{\oplus }\limits_{n=1}^{\infty } W_{2}^{1} \left(H_{n} ,\Delta _{n} \right), \mathop{W_{2}^{1} }\limits^{0} =\mathop{\oplus }\limits_{n=1}^{\infty } \mathop{W_{2}^{1} }\limits^{0} \left(H_{n} ,\Delta _{n} \right),{\rm \; }H=\mathop{\oplus }\limits_{n=1}^{\infty } H_{n} ,{\rm \; }cl\left(T\right)$-closure of the operator $T$. $l\left(\cdot \right)$ is a linear multipoint differential-operator expression for first order in $L^{2} $ in the following form

\begin{equation} \label{GrindEQ__2_1_}
l\left(u\right)=\left(l_{n} \left(u_{n} \right)\right)
\end{equation}
and for each $n\ge 1$
\begin{equation} \label{GrindEQ__2_2_}
l_{n} \left(u_{n} \right)=u'_{n} + A_{n} u_{n} ,
\end{equation}
where $A_{n} :D\left(A_{n} \right)\subset H_{n} \to H_{n} $ is a linear positive defined selfadjoint operator in $H_{n} $.\\
\indent It is clear that formally adjoint expression to (\ref{GrindEQ__2_2_}) in the Hilbert space $L_{n}^{2} $ is in the form\\
\begin{equation} \label{GrindEQ__2_3_}
l_{n}^{+} \left(v_{n} \right)=-v'_{n} + A_{n} v_{n} ,{\rm \; }n\ge 1.
\end{equation}
We define an operator $L'_{n0} $ on the dense manifold of vector functions $D'_{n0} $ in $L_{n}^{2} $ as
\[D^{'}_{n0} :=\left\{{u_{n} \in L_{n}^{2}} :u_{n} =\sum _{k=1}^{m}\phi _{k} f_{k} ,{\rm \; }\phi _{k} \in C_{0}^{\infty } ({\Delta _{n}}),{\rm \;} f_{k} \in D(A_{n}),k=1,2,\cdots ,m;{\rm \; }m\in {\mathbb N} \right\}\]
\noindent with $L'_{n0} u_{n} :=l_{n} \left(u_{n} \right),{\rm \; }n\ge 1$.\\
Since the operator $A_{n}>0, n\geq1$, then from the relation\\
\[Re\left(L'_{n0} u_{n} ,u_{n} \right)_{L_{n}^{2} } =2\left(A_{n} u_{n} ,u_{n} \right)_{L_{n}^{2} } \ge 0,{\rm \; }u_{n} \in D'_{n0} \]
it implies that $L'_{n0} $ is an accretive in $L_{n}^{2} ,{\rm \; }n\ge 1$. Hence the operator $L'_{n0} $ has a closure in $L_{n}^{2} ,{\rm \; }n\ge 1$. The closure $cl\left(L'_{n0} \right)$ of the operator $L'_{n0} $ is called the minimal operator generated by differential-operator expression (\ref{GrindEQ__2_2_}) and is denoted by $L_{n0} $ in $L_{n}^{2} ,{\rm \; }n\ge 1$. The operator $L_{0}$ is defined by
\[D\left(L_{0} \right):=\left\{u=\left(u_{n}\right):u_{n} \in D\left(L_{n0} \right),{\rm \; }n\ge 1,{\rm \; }\sum _{n=1}^{\infty }\left\| L_{n0} u_{n} \right\| _{L_{n}^{2} }^{2} <+\infty  \right\}\]
\noindent with
\[L_{0} u:=\left(L_{n0}u_{n}\right),{\rm \; }u\in D\left(L_{0} \right), L_{0} :D\left(L_{0} \right)\subset L^{2} \to L^{2} \]
is called a minimal operator ( multipoint ) generated by differential-operator expression (\ref{GrindEQ__2_1_}) in Hilbert space $L^{2} $ and denoted by $L_{0} =\mathop{\oplus }\limits_{n=1}^{\infty } L_{n0} $.\\
\indent In a similar way the minimal operator for twopoints  denoted by $L_{n0}^{+} $ in $L_{n}^{2} ,{\rm \; }n\ge 1$ for the formally adjoint linear differential-operator expression (\ref{GrindEQ__2_3_}) can be constructed.\\
\indent In this case the operator  $L_{0}^{+} $ defined by
\[D\left(L_{0}^{+} \right):=\left\{v:=\left(v_{n} \right):v_{n} \in D\left(L_{n0}^{+} \right),{\rm \; }n\ge 1,{\rm \; }\sum _{n=1}^{\infty }\left\| L_{n0}^{+} v_{n} \right\| _{L_{n}^{2} }^{2} <+\infty  \right\}\]
\noindent with $L_{0}^{+} v:=\left(L_{n0}^{+} v_{n} \right),{\rm \; }v\in D\left(L_{0}^{+} \right), L_{0}^{+} :D\left(L_{0}^{+} \right)\subset L^{2} \to L^{2}$
is called a minimal operator (multipoint) generated by $l^{+} \left(v\right)=(l_{n}^{+}(v_{n}))$ in the Hilbert space $L^{2} $ and denoted by $L_{0}^{+} =\mathop{\oplus }\limits_{n=1}^{\infty } L_{n0}^{+}$.\\
\indent We now state the following relevant result.\\
\noindent \textbf{Theorem 2.1.} The minimal operators $L_{0} $ and $L_{0}^{+} $ are densely defined closed operators in $L^{2} $.\\
\indent The following defined operators in $L^{2} $
  $L:=\left(L_{0}^{+} \right)^{*} =\mathop{\oplus }\limits_{n=1}^{\infty } L_{n} $ and $L^{+} :=\left(L_{0} \right)^{*} =\mathop{\oplus }\limits_{n=1}^{\infty } L_{n}^{+} $\\
\noindent are called maximal operators (multipoint) for the differential-operator expression $l\left(\cdot \right)$ and $l^{+} \left(\cdot \right)$ respectively.It is clear that $Lu=\left(l_{n} \left(u_{n} \right)\right),{\rm \; }u\in D\left(L\right)$,
\[D\left(L\right):=\left\{u=\left(u_{n} \right)\in L^{2} :u_{n} \in D\left(L_{n} \right),{\rm \; }n\ge 1{\rm \; ,\; }\sum _{n=1}^{\infty }\left\| L_{n} u_{n} \right\| _{L_{n}^{2} }^{2} <\infty  \right\},\]
\[L^{+} v=\left(l_{n}^{+} \left(v_{n} \right)\right),{\rm \; }v\in D\left(L^{+} \right),\]
\[D\left(L^{+} \right):=\left\{v=\left(v_{n} \right)\in L^{2} :v_{n} \in D\left(L_{n}^{+} \right),{\rm \; }n\ge 1{\rm \; ,\; }\sum _{n=1}^{\infty }\left\| L_{n}^{+} v_{n} \right\| _{L_{n}^{2} }^{2} <\infty  \right\}\]
and $L_{0} \subset L$, $L_{0}^{+} \subset L^{+} $.\\
\indent Furthermore, the validity of following proposition is clear.

\noindent \textbf{Theorem 2.2.} The domain of the operators $L$ and $L_{0} $ are
\[\begin{array}{l} {D\left(L\right)=\left\{u=\left(u_{n} \right)\in L^{2} :\left(1\right){\rm \; for\; each\; }n\ge 1{\rm \; vector\; function\; }u_{n} \in L_{n}^{2} ,{\rm \; \; } u_{n} \right. } \\ {{\rm \; \; \; \; \; \; \; \; \; \; \; \; \; \; \; \; \; \; \; \; \; \; \; \; \; \; \; \; \; \; \; \; \; \; \; \; \; \; \; \; is\; absolutely\; continuous\; in\; interval\; }\Delta _{n} ;} \\ {{\rm \; \; \; \; \; \; \; \; \; \; \; \; \; \; \; \; \; \; \; \; \; \;\; \; \; \; \; \; \; \; \; \; \;\; \; \; \; \; \; }\left(2\right){\rm \; }l_{n} \left(u_{n} \right)\in L_{n}^{2} ,{\rm \; }n\ge 1;{\rm \; }\left. \left(3\right){\rm \; }l\left(u\right)=\left(l_{n} \left(u_{n} \right)\right)\in L^{2} \right\}} \\
{{\rm \; \; \; \; \; \; \; \; \; }=\left\{u=\left(u_{n} \right)\in L^{2} :u_{n} \in D\left(L_{n} \right),{\rm \; }n\ge 1{\rm \; and\; }l\left(u\right)=\left(l_{n} \left(u_{n} \right)\right)\in L^{2} \right\},} \end{array}\]
$ D(L_{0})=\left\{u=\left(u_{n} \right)\in D\left(L\right):u_{n} \left(a_{n} \right)=u_{n} \left(b_{n} \right)=0,{\rm \; }n\ge 1\right\}$.\\
\noindent \textbf{Remark 2.3.} If $A_{n} \in B\left(H\right),{\rm \; }n\ge 1$ and $\mathop{\sup }\limits_{n\ge 1} \left\| A_{n} \right\| \le c<+\infty $, then for any $u=\left(u_{n} \right)\in L^{2} $ we have $\left(Au\right)=\left(A_{n} u_{n} \right)\in L^{2} $.

Now the following results can be proved .

\noindent \textbf{Theorem 2.4.} If a minimal operator $L_{0} $ is formally normal in $L^{2} $, then $D\left(L_{0} \right)\subset \mathop{W_{2}^{1} }\limits^{0} $ and  $AD\left(L_{0} \right)\subset L^{2} $.

\noindent \textbf{Theorem 2.5.} If $A^{1/2} W_{2}^{1} \subset W_{2}^{1} $, then minimal operator $L_{0} $ is formally normal in $L^{2} $.

\noindent \textbf{Proof:} In this case from the following relations
\[L_{0}^{+} u=L_{0} u-2Au,{\rm \; }u\in D\left(L_{0} \right),\]
\[L_{0} u=L_{0}^{+} u+2Au,{\rm \; }u\in D\left(L_{0}^{+} \right)\]
it implies that $D\left(L_{0} \right)=D\left(L_{0}^{+} \right)$. Since $D\left(L_{0}^{+} \right)\subset D\left(L_{0}^{*} \right)=D\left(L^{+} \right)$, it is obtained that $D\left(L_{0} \right)\subset D\left(L^{+} \right)$.\\
 On the other hand for any $u\in D\left(L_{0} \right)$
\[\begin{array}{l} {\left\| L_{0} u\right\| _{L^{2} }^{2} =\left(u'+Au,u'+Au\right)_{L^{2} } =\left\| u'\right\| _{L^{2} }^{2} +\left[\left(u',Au\right)_{L^{2} } +\left(Au,u'\right)_{L^{2} } \right]+\left\| Au\right\| _{L^{2} }^{2} } \\ {{\rm \; \; \; \; \; \; \; \; \; \; \; }=\left\| u'\right\| _{L^{2} }^{2} +\left\| Au\right\| _{L^{2} }^{2} } \end{array}\]
and
\[\begin{array}{l} {\left\| L^{+} u\right\| _{L^{2} }^{2} =\left(-u'+Au,-u'+Au\right)_{L^{2} } =\left\| u'\right\| _{L^{2} }^{2} -\left[\left(u',Au\right)_{L^{2} } +\left(Au,u'\right)_{L^{2} } \right]+\left\| Au\right\| _{L^{2} }^{2} } \\ {{\rm \; \; \; \; \; \; \; \; \; \; \; }=\left\| u'\right\| _{L^{2} }^{2} +\left\| Au\right\| _{L^{2} }^{2} .} \end{array}\]
Thus, it is established that operator $L_{0} $ is formally normal in $L^{2} $.

\noindent \textbf{Remark 2.6.} If  $A_{n} \in B\left(H\right),{\rm \; }n\ge 1$ and $\mathop{\sup }\limits_{n\ge 1} \left\| A_{n} \right\| \le c<+\infty $, then $D\left(L_{0} \right)=D\left(L_{0}^{+} \right)$ and $D\left(L\right)=D\left(L^{+} \right).$

\noindent \textbf{Remark 2.7.} If $AW_{2}^{1} \subset L^{2} $, then $D\left(L_{0} \right)=D\left(L_{0}^{+} \right)$ and $D\left(L\right)=D\left(L^{+} \right)$.

\section{\textbf{Description of Normal Extensions of the Minimal Operator}}

 In this section the main purpose is to describe all normal extensions of the minimal operator  \textbf{$L_{0} $ }in\textbf{ $L^{2} $ }in terms in the boundary values of the endpoints of the\textbf{ }subintervals .

  First, we will show that there exists normal extension of the minimal operator \textbf{$L_{0} $}. Consider the following extension of the minimal operator $L_{0} $
\[\left\{\begin{array}{l} {{\rm \; \; \; \; }\widetilde{L}u:=u'+Au,{\rm \; }AW_{2}^{1} \subset W_{2}^{1} ,} \\ {D\left(\widetilde{L}\right)=\left\{u=\left(u_{n} \right)\in W_{2}^{1} :u_{n} \left(a_{n} \right)=u_{n} \left(b_{n} \right),{\rm \; }n\ge1\right\}.} \end{array}\right. \]
 Under the condition on the coefficient $A$ we have
\[\begin{array}{l} {{\left(\widetilde{L}u,v\right)_{L^{2} } =\left(u',v\right)_{L^{2} } +\left(Au,v\right)_{L^{2} }}}=\left(u,v\right)_{L^{2}} ^{{'} } +\left(u,-v^{'}+Av\right)_{L^{2}}\\
{{\rm \; \; \; \; \; \; \; }=\sum \limits_{n=1}^{\infty }[(u_{n}(b_{n}),v_{n}(b_{n} ))_{H_{n}}- (u_{n}(a_{n}),v_{n}(a_{n} ))_{H_{n}}]+(u,-v^{'}+Av)_{L^{2}}}\end{array}\]
From this it is obtained \\
\[\left\{\begin{array}{l} {{\rm \; \; \; \; }\widetilde{L}^{*}}v:=-v'+Av,\\
{D\left(\widetilde{L}^{*}\right)=\left\{v=\left(v_{n} \right)\in W_{2}^{1} :v_{n} \left(a_{n} \right)=v_{n} \left(b_{n} \right),{\rm \; }n\ge 1\right\}.} \end{array}\right. \]
\noindent In this case it is clear that $D\left(\widetilde{L}\right)=D\left(\widetilde{L}^{*} \right)$. On the other hand, since for each $u\in D\left(\widetilde{L}\right)$
\[\begin{array}{l} {\left\| \widetilde{L}u\right\| _{L^{2} }^{2} =\left\| u'\right\| _{L^{2} }^{2} +\left[\left(u',Au\right)_{L^{2} } +\left(Au,u'\right)_{L^{2} } \right]+\left\| Au\right\| _{L^{2} }^{2} ,} \\ {\left\| \widetilde{L}^{*} u\right\| _{L^{2} }^{2} =\left\| u'\right\| _{L^{2} }^{2} -\left[\left(u',Au\right)_{L^{2} } +\left(Au,u'\right)_{L^{2} } \right]+\left\| Au\right\| _{L^{2} }^{2} } \end{array}\]
and \\
\[\begin{array}{l} {\left(u',Au\right)_{L^{2} } +\left(Au,u'\right)_{L^{2} } =\left(u,Au\right)_{L^{2} } ^{{'} } } \\ \noindent{=\sum\limits_{n=1}^{\infty}[(u_{n}(b_{n}),A_{n}u_{n}(b_{n}))_{H_{n}}-(u_{n}(a_{n}),A_{n}u_{n}(a_{n}))_{H_{n}} ]=0}.\\
\end{array}\]
Then $\left\| \widetilde{L}u\right\| _{L^{2} } =\left\| \widetilde{L}^{*} u\right\| _{L^{2} } $ for every $u\in D\left(\widetilde{L}\right)$. Consequently, $\widetilde{L}$ is a normal extension of the  minimal operator $L_{0} $.

 The following result establishes the relationship between normal extensions of $L_{0} $ and normal extensions of  $L_{n0} ,{\rm \; }n\ge 1$.

\noindent \textbf{Theorem 3.1.} The extension $\widetilde{L}=\mathop{\oplus }\limits_{n=1}^{\infty }\widetilde{L_{n}}$ of the minimal operator $L_{0}$ in $ L^{2}$ is a normal if and only if for any $n\geq 1$, $\widetilde{L_{n}}$ is so in $L_{n}^{2}$ .

 Now using the Theorem 3.1 and \cite{Ism2} we can formulate the following main result of this section, where it is given a description of all normal extension  of the minimal operator $L_{0} $ in $L^{2} $ in terms of boundary values of vector functions at the endpoints of subintervals.\\
\noindent \textbf{Theorem 3.2.} Let $A^{1/2} W_{2}^{1} \subset W_{2}^{1} $. If $\widetilde{L}=\mathop{\oplus }\limits_{n=1}^{\infty }\widetilde{L_{n} }$ is a normal extension of the minimal operator $L_{0} $ in $L^{2} $, then it is generated by differential-operator expression (\ref{GrindEQ__2_1_}) with boundary conditions
\begin{equation} \label{GrindEQ__3_0_}
u_{n}(b_{n})=W_{n}u_{n}(a_{n}), u_{n}\in D(L_{n}),
\end{equation}
where $W_{n}$ is a unitary operator in $H_{n}$ and  $W_{n}A_{n}^{-1}=A_{n}^{-1}W_{n}, n\geq 1$. The unitary operator $W=\mathop{\oplus }\limits_{n=1}^{\infty }W_{n}$ in $H=\mathop{\oplus }\limits_{n=1}^{\infty }H_{n}$ is determined uniquely by the extension $\widetilde{L}$, that is $\widetilde{L}=L_{W} $.\\
\indent On the contrary, the restriction of the maximal operator $L$ to the linear manifold $u\in D(L)$ satisfying the condition (\ref{GrindEQ__3_0_}) with any unitary operator $W=\mathop{\oplus }\limits_{n=1}^{\infty }W_{n}$ in $H$  with property $WA^{-1}=A^{-1}W$ is a normal extension of the minimal operator $L_{0}$ in $L^{2}$.

\section{\textbf{Some Compactness Properties of The Normal Extensions}}

 The following two proposition can be easily proved in general case.

\noindent \textbf{Theorem 4.1.} For the point spectrum of ${\mathscr{A}}=\mathop{\oplus }\limits_{n=1}^{\infty }\mathscr{A}_{n}$ in the direct sum ${\mathscr{H}}=\mathop{\oplus }\limits_{n=1}^{\infty }\mathscr{H}_{n}$ of Hilbert spaces $\mathscr{H}_{n}, n\geq1$ ıt is true that
\[\sigma_{p}(\mathscr{A})=\bigcup\limits_{n=1}^{\infty}\sigma_{p}(\mathscr{A}_{n})\]

\noindent \textbf{Theorem 4.2.} Let $\mathscr{A}_{n}\in B(\mathscr{H}_{n}), n\geq 1, \mathscr{A}=\mathop{\oplus }\limits_{n=1}^{\infty }\mathscr{A}_{n}$ and
$\mathscr{H}=\mathop{\oplus }\limits_{n=1}^{\infty }\mathscr{H}_{n}$. In order for $\mathscr{A}\in B(\mathscr{H})$ the necessary and sufficient condition is that the $\underset{{n\geq1}}{sup}\|\mathscr{A}_{n}\|$ be finite. In this case $\|\mathscr{A}\|= \underset{{n\geq1}}{sup}\|\mathscr{A}_{n}\|$.\\

\indent Let $C_{\infty}(\cdot)$ and $C_{p}(\cdot), 1\leq p<\infty$ denote the class of compact operators and the Schatten-von Neumann subclasses of compact operators in corresponding spaces respectively.

\noindent \textbf{Definition 4.3.}\cite{Gor} Let $T$ be a linear closed and densely defined operator in any Hilbert space $\mathfrak{H}$. If $\rho(T)\neq{\O}$ and for $\lambda\in\rho(T)$ the resolvent operator $R_{\lambda}(T)\in C_{\infty}(\mathfrak{H})$, then operator $T:D(T)\subset{\mathfrak{H}}\to{\mathfrak{H}}$ is called an operator with discrete spectrum  \\

\indent In first note that the following results are true.\\
\noindent \textbf{Theorem 4.4.} If the operator $\mathscr{A}=\mathop{\oplus }\limits_{n=1}^{\infty }\mathscr{A}_{n}$ as an operator with discrete spectrum in $\mathscr{H}=\mathop{\oplus }\limits_{n=1}^{\infty }\mathscr{H}_{n}$, then for every $n\geq1$ the operator $\mathscr{A}_{n}$ is so in $\mathscr{H}_{n}$.\\

\noindent \textbf{Remark 4.5.} Unfortunately, the converse of the Theorem 4.4 is not true in general case.\\
\noindent Indeed, consider the following sequence of operators $\mathscr{A}_{n}u_{n}=u_{n}, 0<dim\mathscr{H}_{n}=d_{n}<\infty, n\geq1$. In this case for every $n\geq1$ operator $\mathscr{A}_{n}$ is an operator with discrete spectrum. But an inverse of the direct sum operator $\mathscr{A}=\mathop{\oplus }\limits_{n=1}^{\infty }\mathscr{A}_{n}$ is not compact operator in $\mathscr{H}=\mathop{\oplus }\limits_{n=1}^{\infty }\mathscr{H}_{n}$, because $dim\mathscr{H}=\infty$ and $\mathscr{A}$ is an identity operator in $\mathscr{H}$.\\

\noindent \textbf{Theorem 4.6.} If $\mathscr{A}=\mathop{\oplus }\limits_{n=1}^{\infty }\mathscr{A}_{n}, \mathscr{A}_{n}$ is an operator with discrete spectrum in $\mathscr{H}_{n}$, $n\geq1$, $\bigcap\limits_{n=1}^{\infty}{\rho(\mathscr{A}_{n})}\neq{\O}$ and $\lim\limits_{n\to\infty}\|R_{\lambda}(\mathscr{A}_{n})\|=0$, then $\mathscr{A}$ is an operator with discrete spectrum in $\mathscr{H}$.\\

\noindent \textbf{Proof.} For each $\lambda \in \bigcap\limits_{n=1}^{\infty}{\rho(\mathscr{A}_{n})}$ we have $R_{\lambda}(\mathscr{A}_{n})\in C_{\infty}(\mathscr{H}_{n}), n\geq1$.\\
\indent Now define the following operators $\mathscr{K}_{m}:\mathscr{H}\to\mathscr{H}, m\geq1$ as
\[\mathscr{K}_{m}:=\{R_{\lambda}(\mathscr{A}_{1})u_{1}, R_{\lambda}(\mathscr{A}_{2})u_{2},\dots, R_{\lambda}(\mathscr{A}_{m})u_{m}, 0, 0, \dots \}, u=(u_{n})\in \mathscr{H}. \]
\indent The convergence of the operators $\mathscr{K}_{m}$ to the operator $\mathscr{K}$ in operator norm will be investigated. For the $u=(u_{n})\in \mathscr{H}$ we have
\begin{eqnarray*}
\|\mathscr{K}_{m}u-\mathscr{K}u\|_{\mathscr{H}}^{2} & = & \sum_{n=m+1}^{\infty}\|R_{\lambda}(\mathscr{A}_{n})u_{n}\|_{\mathscr{H}_{n}}^{2}
\leq\sum\limits_{n=m+1}^{\infty}\|R_{\lambda}(\mathscr{A}_{n})\|^{2}\|\|u_{n}\|_{\mathscr{H}_{n}}^{2}\\
& \leq & \left(\underset{n\geq m+1}{sup}{\|R_{\lambda}(\mathscr{A}_{n})\|}\right)^{2}\sum_{n=1}^{\infty}\|u_{n}\|_{\mathscr{H}_{n}}^{2}=\left(\underset{n\geq m+1}{sup}{\|R_{\lambda}(\mathscr{A}_{n})\|} \right)^{2}\|u\|_{\mathscr{H}}^{2}
\end{eqnarray*}
\noindent thus we get $\|\mathscr{K}_{m}u-\mathscr{K}u\|\leq \underset{n\geq m+1}{sup}{\|R_{\lambda}(\mathscr{A}_{n})\|}, m\geq1$.\\
\indent This means that sequence of operators $(\mathscr{K}_{m})$ converges in operator norm to the operator $\mathscr{K}$. Then by the important theorem of the theory of compact operators it is implies that $\mathscr{K}\in C_{\infty}(\mathscr{H})$ \cite{Dun}, because for any $m\geq1$, $\mathscr{K}_{m}\in C_{\infty}(\mathscr{H})$.\\
\indent Finally, using the Theorem 4.6 can be proved the following result.\\

\noindent \textbf{Theorem 4.7.} If $A_{n}^{-1}\in C_{\infty}(\mathscr{H}_{n}), n\geq1$, $\underset{n\geq1}{sup}(b_{n}-a_{n})<\infty$ and the sequence of first minimal eigenvalues $\lambda_{1}(\mathscr{A}_{n})$ of the operators $\mathscr{A}_{n}, n\geq1$ is satisfy the condition
\[\lambda_{1}(\mathscr{A}_{n})\to\infty {\rm \;\;\;as}{\rm \;\;\;} n \to \infty,  \]
\noindent then the extension $\widetilde{L}=\mathop{\oplus }\limits_{n=1}^{\infty }{L_{n}}$ is an operator with discrete spectrum in $L^{2}$.\\

\noindent \textbf{Theorem 4.8.} Let $H=\mathop{\oplus }\limits_{n=1}^{\infty }{H_{n}}, A=\mathop{\oplus }\limits_{n=1}^{\infty }{A_{n}}$ and $A_{n}\in C_{p}(H_{n}), n\geq1, 1\leq p<\infty$. In order for $A\in C_{p}(H)$ the necessary and sufficient condition is that the series $\sum\limits_{n=1}^{\infty }{\sum\limits_{k=1}^{\infty }}\mu_{k}^{p}(A_{n})$ be convergent.\\
\indent Now we will dedicate an application the last theorem.\\
\indent For all $n\geq1, \mathfrak{H}_{n}$ is a Hilbert space, $\Delta_{n}=(a_{n},b_{n}), -\infty<a_{n}<b_{n}<a_{n+1}<\cdots<\infty, A_{n}:D(A_{n})\subset \mathfrak{H}_{n}\to \mathfrak{H}_{n}, A_{n}=A_{n}^{*}\geq E, W_{n}:\mathfrak{H}_{n}\to \mathfrak{H}_{n}$ is unitary operator, $A_{n}^{-1}W_{n}=W_{n}A_{n}^{-1}, L_{W_n}u_{n}=u'_{n}+A_{n}u_{n}, A_{n}W_{2}^{1}(\mathfrak{H}_{n},\Delta _{n})\subset W_{2}^{1}(\mathfrak{H}_{n},\Delta _{n}), H_{n}=L^{2}(\mathfrak{H}_{n},\Delta _{n}), D(L_{W_n})=\{u_{n}\in W_{2}^{1}(\mathfrak{H}_{n},\Delta _{n}): u_{n}(b_{n})=W_{n}u_{n}(a_{n})\}, L_{W_n}:H_{n}\to H_{n}, W=\mathop{\oplus }\limits_{n=1}^{\infty }{W_{n}}, L_{W}=\mathop{\oplus }\limits_{n=1}^{\infty }{L_{W_n}}, H=\mathop{\oplus }\limits_{n=1}^{\infty }{H_{n}}$ and $h=\underset{n\geq 1}{sup}(b_{n}-a_{n})<\infty$.\\
\indent Since for all $n\geq1${\rm \;\;} $W_{n}$ is a unitary operator in $\mathfrak{H}_{n}$, then $L_{W_n}$ is normal operator in $H_{n}$ \cite{Ism2}. Also for $L_{W}:D(L_{W})\subset H\to H$,{\rm \;\;}the relation $L_{W}L_{W}^{*}=L_{W}^{*}L_{W}$ is true, i.e. $L_{W}$ is a normal operator in $H$. It is known that, if $A_{n}^{-1}\in C_{p}(\mathfrak{H}_{n})$, for $p>1$, then $L_{W_n}^{-1}\in C_{p}(H_{n}), p>1$ for all $n\geq1$\cite{Ism2}. On the other hand, if $A_{n}^{-1}\in C_{\infty}(\mathfrak{H}_{n}), n\geq1$, then eigenvalues $\lambda_{q}(L_{W_n}), q\geq1$ of operator $L_{W_n}$ is in form\\
\[\lambda_{q}(L_{W_n})=\lambda_{m}(A_{n})+\frac{i}{a_{n}-b_{n}}\left( arg\lambda_{m}(W_{n}^{*}e^{(-A_{n}(b_{n}-a_{n}))})+2k\pi \right), m\geq1, k\in{\mathbb Z}, n\geq1 ,\]
\noindent where $q=q(m,k)\in {\mathbb N}, m\geq1, k\in {\mathbb Z}$. Therefore we have the following corollary.\\

\noindent \textbf{Theorem 4.9.} If $A=\mathop{\oplus }\limits_{n=1}^{\infty}A_{n}, \mathfrak{H}=\mathop{\oplus }\limits_{n=1}^{\infty}\mathfrak{H}_{n}$ and $A^{-1}\in C_{p/2}(\mathfrak{H}), 2<p<\infty$, then $L_{W}^{-1}\in C_{p}(H)$.\\

\noindent \textbf{Proof.} The operator $L_{W}$ is a normal in $H$. Consequently, for the characteristic numbers of normal operator $L_{W}^{-1}$ an equality $\mu_{q}(L_{W}^{-1})=\mid\lambda_{q}(L_{W}^{-1})\mid, q\geq1$ holds \cite{Dun}. Now we search for convergence of the series $\sum_{q=1}^{\infty}\mu_{q}^{p}(L_{W}^{-1}), 2<p<\infty$.
\begin{eqnarray*}
 \sum\limits_{n=1}^{\infty }\sum\limits_{q=1}^{\infty }\mu_{q}^{p}(L_{W_n}^{-1}) & = & \sum\limits_{n=1}^{\infty }\sum\limits_{k=-\infty}^{\infty}\sum\limits_{m=1}^{\infty }\left( \lambda_{m}^{2}(A_{n})+\frac{1}{(b_{n}-a_{n})^{2}}(\delta(m,n)+2k\pi)^{2} \right)^{-p/2} \\
& \leq & \sum\limits_{n=1}^{\infty }\sum\limits_{k=-\infty}^{\infty}\sum\limits_{m=1}^{\infty }\left( \lambda_{m}^{2}(A_{n})+\frac{4k^{2}\pi^{2}}{(b_{n}-a_{n})^{2}} \right)^{-p/2} \\
& \leq & \sum\limits_{n=1}^{\infty }\sum\limits_{m=1}^{\infty}\left( \lambda_{m}^{2}(A_{n}) \right)^{-p/2}+2\sum\limits_{n=1}^{\infty }\sum\limits_{k=1}^{\infty}\sum\limits_{m=1}^{\infty}\left( \lambda_{m}^{2}(A_{n})+\frac{4k^{2}\pi^{2}}{(b_{n}-a_{n})^{2}} \right)^{-p/2}
\end{eqnarray*}
\noindent where $\delta(m,n)=arg\lambda_{m}(W_{n}^{*}e^{(-A_{n}(b_{n}-a_{n}))}), n\geq1, m\geq1$. Then from the inequality $\frac{\mid ts \mid}{t^{2}+s^{2}}\leq \frac{1}{2}$ for all $t,s\in \mathbb{R}\setminus \{0\}$ and last equation we have the inequality
\[\sum\limits_{n=1}^{\infty }\sum\limits_{k=1}^{\infty}\sum\limits_{m=1}^{\infty}\left( \lambda_{m}^{2}(A_{n})+\frac{4k^{2}\pi^{2}}{(b_{n}-a_{n})^{2}} \right)^{-p/2}\leq 2^{-p}\pi^{-p/2}h^{p/2}\left( \sum\limits_{n=1}^{\infty }\sum\limits_{m=1}^{\infty}\left\vert\frac{1}{\lambda_{m}(A_{n})}\right\vert^{p/2}\sum\limits_{k=1}^{\infty}\left\vert \frac{1}{k} \right\vert^{p/2}\right)\]
\noindent Since $A^{-1}\in C_{p/2}(\mathfrak{H})$, then the series $\sum\limits_{n=1}^{\infty }\sum\limits_{m=1}^{\infty}\left\vert \lambda_{m}(A_{n}) \right\vert^{-p/2}$ is convergent. Thus the series
\[\sum\limits_{n=1}^{\infty }\sum\limits_{k=1}^{\infty}\sum\limits_{m=1}^{\infty}\left( \lambda_{m}^{2}(A_{n})+\frac{4k^{2}\pi^{2}}{(b_{n}-a_{n})^{2}} \right)^{-p/2} \]
is also convergent. Then from the relation
\[\sum\limits_{n=1}^{\infty }\sum\limits_{m=1}^{\infty}\left\vert \lambda_{m}(A_{n})\right\vert^{-p}\leq\sum\limits_{n=1}^{\infty }\sum\limits_{m=1}^{\infty}\left\vert \lambda_{m}(A_{n})\right\vert^{-p/2}\]
\noindent and the convergence of the series $\sum\limits_{n=1}^{\infty }\sum\limits_{m=1}^{\infty}\left\vert \lambda_{m}(A_{n})\right\vert^{-p/2}$ we get that the series $\sum\limits_{n=1}^{\infty }\sum\limits_{m=1}^{\infty}\left\vert \lambda_{m}(A_{n})\right\vert^{-p}$ is convergent too. Consequently the series $\sum\limits_{q=1}^{\infty}\mu_{q}^{p}(L_{W}^{-1}), 2<p<\infty$ is convergent and thus $L_{W}^{-1}\in C_{p}(H), 2<p<\infty$.\\
\indent The Theorem 4.8 and 4.9. can be can generalized.  \\
\noindent \textbf{Corollary 4.10.} Let for all $n\geq1${\rm \:}$A_{n}\in C_{p_n}(H_{n}), 1\leq p_{n}<\infty$ and $p=\underset{n\geq1}{sup}{\rm\;}p_{n}<\infty$. For $A\in C_{p}(H)$ the necessary and sufficient condition is that the series $\sum\limits_{n=1}^{\infty }\sum\limits_{k=1}^{\infty}\mu_{k}^{p}(A_{n})$ be convergent.\\
\noindent \textbf{Theorem 4.11.} If $A_{n}^{-1}\in C_{{p_n}/2}(H_{n}), 2\leq p_{n}<\infty, p=\underset{n\geq1}{sup}{\rm\;}p_{n}<\infty$, then $L_{W}^{-1}\in C_{p}(H)$.

\footnotesize
\noindent \textbf{}

\centerline{}\centerline{}
\end{document}